\documentclass[12pt,reqno]{amsart} 
\usepackage[latin1]{inputenc}
\usepackage{amsmath, amsthm, amssymb,color,graphicx,bm,psfrag,float,verbatim}
\numberwithin{equation}{section}
\usepackage[all]{xy} 
\usepackage{fullpage} 
\newcommand{\ra}{\rightarrow} 
\newcommand{\lra}{\longrightarrow} 
\renewcommand{\P}{\mathbb P} 
\newcommand{\conic}{\mathcal K}
\newcommand{\complex}{\mathbb C} 
\newcommand{\Sym}{\text{Sym}}
\newcommand{\InvH}{\mathcal Y} 
\newcommand{\RicoV}{\mathcal R} 
\newcommand{\scD}{\mathcal D}
\newcommand{\six}{\text{\sc six}}
\newcommand{\ltr}{\text{\sc ltr}}
\newcommand{\bA}{\mathbb A} 
\newcommand{\bB}{\mathbb B} 
\newcommand{\bC}{\mathbb C} 
\newcommand{\bD}{\mathbb D} 
\newcommand{\bE}{\mathbb E} 
\newcommand{\bF}{\mathbb F} 
\newcommand{\SG}{\mathfrak S} 
\newcommand{\ux}{\mathbf x} 
\newcommand{\germf}{\mathfrak f}
 
\newcommand{\T}{\mathbb T} 
\newcommand{\sz}{\mathsf z}
\newcommand{\Hex}{\mathsf {HEX}}
\newcommand{\Label}{\mathcal L}
\renewcommand{\ge}{\geqslant} 
\renewcommand{\le}{\leqslant} 
\renewcommand{\proof}{{\sc Proof. \;}}

\newtheorem{Theorem}{Theorem}[section]

\newtheorem{Proposition}[Theorem]{Proposition}

\newcommand{\arr}[6]{\left[ \begin{array}{ccc} #1 & #2 & #3\\ #4 & #5 & #6 \end{array} \right]}
\newcommand{\pasc}[6]{\left\{ \begin{array}{ccc} #1 & #2 & #3\\ #4 &
     #5 & #6 \end{array} \right\}}
\newcommand{\intarr}[4]{\left[ \begin{array}{cc} #1 & #2 \\ #3 &
      #4 \end{array} \right]}

\begin{document} 
\title{On the coincidence of Pascal lines} 
\author{Jaydeep Chipalkatti} 
\maketitle

\bigskip 

\parbox{17cm}{ \small
{\sc Abstract:} Let $\conic$ denote a smooth conic in the complex
projective plane. Pascal's theorem says that, given six points
$A,B,C,D,E,F$ on $\conic$, the three intersection points $AE \cap BF,
AD \cap CF, BD \cap CE$ are collinear. This defines the Pascal line
of the array $\arr{A}{B}{C}{F}{E}{D}$,  and one gets sixty such lines
in general by permuting the points. In this paper we
consider the variety $\Psi$ of sextuples $\{A, \dots, F\}$, for which
some of these Pascal lines coincide. We show that $\Psi$ has two
irreducible components: a five-dimensional component of sextuples in involution,
and a four-dimensional component of the so-called `ricochet
configurations'. This gives a complete synthetic characterisation of
points in $\Psi$. The proof relies upon Gr{\"o}bner basis techniques to
solve multivariate polynomial equations.} 

\bigskip 


\bigskip 


\bigskip 

\tableofcontents

\section{Introduction} 
\thispagestyle{empty} 

\subsection{} 
Fix a smooth conic $\conic$ in the complex projective plane $\P^2$,
and choose six distinct points $A, B, C, D, E, F$ on $\conic$. If they are displayed
as an array $\arr{A}{B}{C}{F}{E}{D}$, then Pascal's theorem says that the three `cross-hair' intersection points 
\[ AE \cap BF, \quad AD \cap CF, \quad BD \cap CE, \] 
are collinear (see Diagram 1). 

\begin{figure}
\includegraphics[width=8cm]{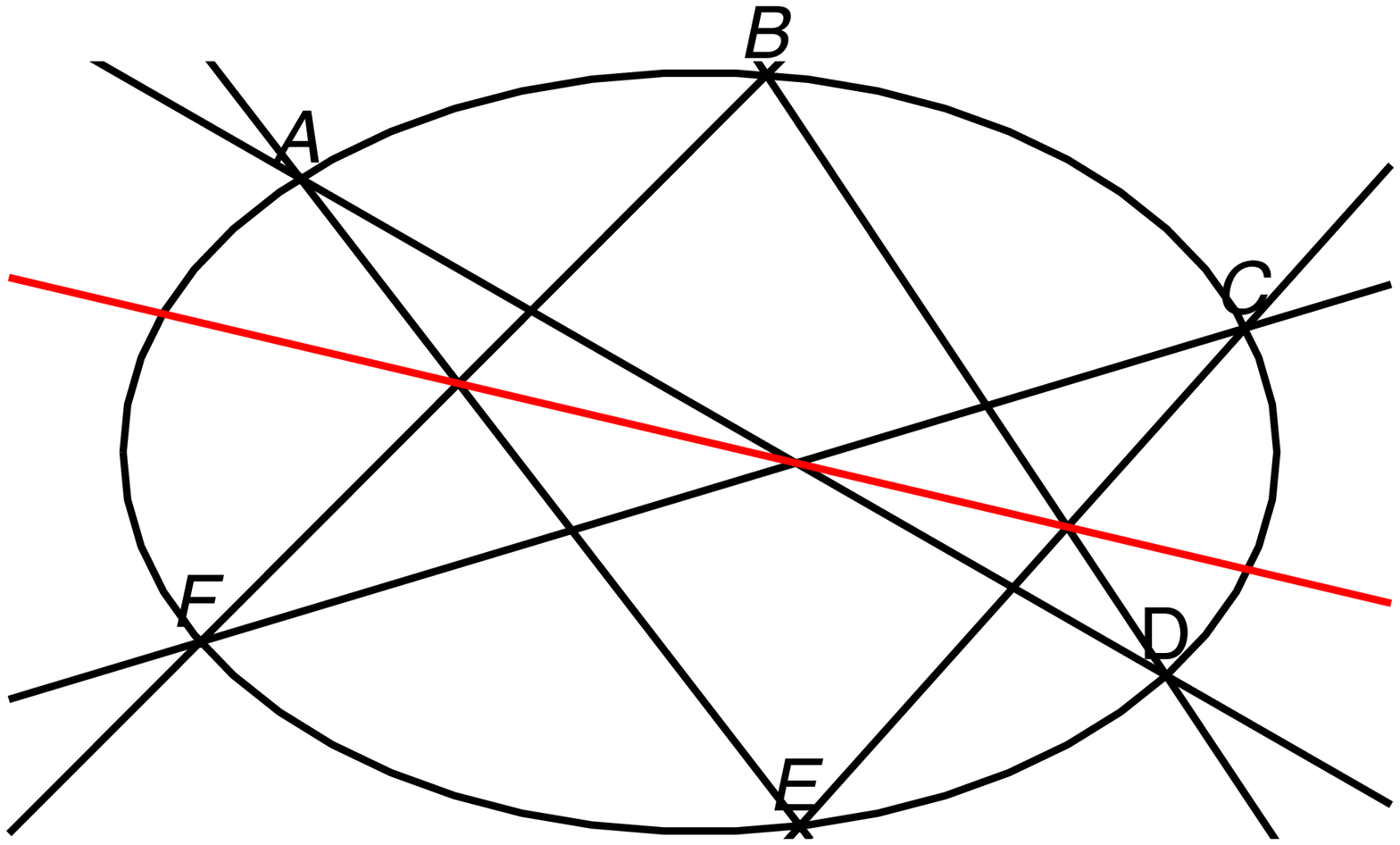}
\centerline{\small \textcolor{magenta}{\bf Diagram 1}}
\end{figure} 

The line containing them (usually called the Pascal line, or just the
Pascal)  will be denoted as $\pasc{A}{B}{C}{F}{E}{D}$. A different
arrangement of the same points, say 
$\pasc{D}{A}{C}{F}{B}{E}$, will~\emph{a priori} give a different
line. A permutation of rows or columns has no effect on intersection points; for instance, 
\[ \pasc{A}{B}{C}{F}{E}{D} = \pasc{F}{E}{D}{A}{B}{C} = \pasc{D}{E}{F}{C}{B}{A} \quad \text{etc.,} \] 
hence one gets at most $6!/(2 \times 3!) = 60$ possibilities for the
Pascal by permuting the points. For a \emph{general} choice of six points, these sixty lines are in fact
distinct (see~\cite{Pedoe}); that is to say, we must be inside a 
special geometric configuration of some kind if any of the Pascals are to coincide. 

\subsection{} 
One such configuration is as follows: suppose that the points are in \emph{involution}, i.e., the
lines $AF, BE, CD$ are concurrent in the point $Q$ (see Diagram 2). 

\begin{figure}
\includegraphics[width=7cm]{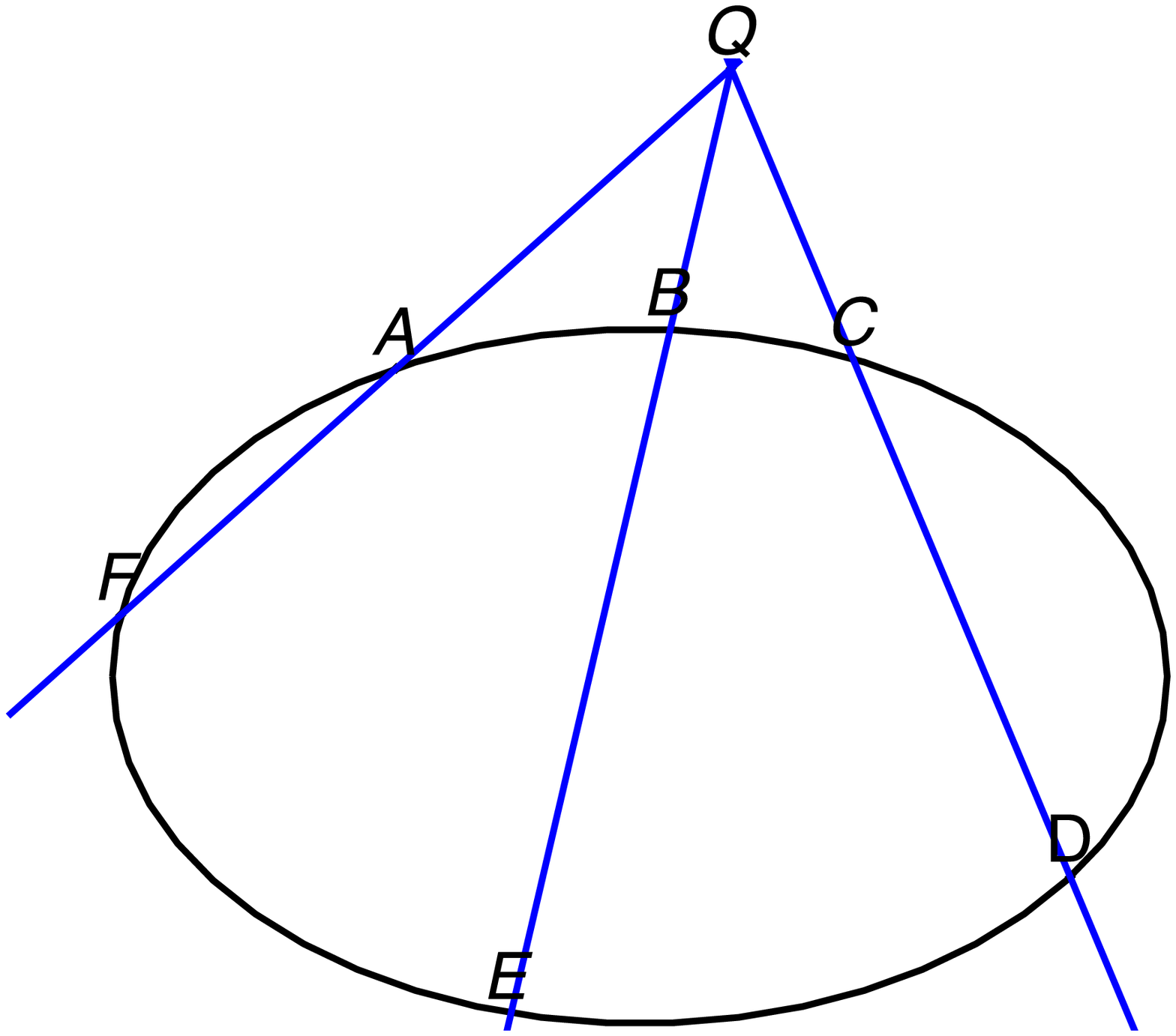}
\centerline{\small \textcolor{magenta}{\bf Diagram 2}}
\end{figure} 

Then it is not difficult to show (see Proposition~\ref{prop.inv.coincidence} below), that the following four
Pascals become equal: 
\begin{equation}
 \pasc{A}{B}{C}{F}{E}{D}, \quad \pasc{A}{B}{D}{F}{E}{C}, \quad 
\pasc{F}{B}{C}{A}{E}{D}, \quad \pasc{A}{E}{C}{F}{B}{D}. 
\label{inv.display} \end{equation}
(The pattern is simple; pick any one column from the first array and
interchange its entries.) 
There are no further coincidences, so that a generic involutive configuration 
has $57$ distinct Pascals. It is natural enough to ask
whether the converse holds, i.e., whether assuming that some two
Pascals coincide forces the initial six points to be in
involution. The main result of this paper (Theorem~\ref{main.theorem}
below) says that the answer is `No, but almost yes.' This requires some
explanation. 

\subsection{} \label{section.pascal.map} 
Since $\conic$ is isomorphic to the projective line $\P^1$, an unordered 
sextuple of points in $\conic$ may be identified with
an element in the symmetric product 
\[ \Sym^6(\P^1) = \frac{(\P^1 \times \P^1 \dots \times
\P^1)}{\text{symmetric group on six objects}} \simeq \P^6. \] 
Let $\Delta \subseteq \P^6$ denote the discriminant hypersurface
parametrising sextuples where the points are not all distinct. Then we have a morphism 
\[ \P^6 \setminus \Delta \stackrel{\germf}{\lra}
\Sym^{60} (\P^2)^*, \] 
which sends a sextuple to all of its Pascals. If $\scD \subseteq \Sym^{60} (\P^2)^*$ denotes the
`big diagonal' parametrising repeated lines, then $\Psi =
\germf^{-1}(\scD)$ is the variety of sextuples of distinct points whose Pascals are not
all distinct. Our main theorem says that $\Psi$ is a union of two irreducible
components $\InvH$ and $\RicoV$, where 
\begin{itemize}
\item 
$\InvH$ is the degree $15$ hypersurface of sextuples in involution, and 
\item 
$\RicoV$ is the four-dimensional variety of sextuples in what will be
called the `ricochet configuration'. 
\end{itemize} 
Since it is $\InvH$ which has the larger dimension, a general sextuple in $\Psi$ is in involution. 

\subsection{} \label{section.ricochet.construction} 
The ricochet configuration (see Diagram 3) has not appeared in literature to
the best of my knowledge. I arrived at it after a measure of 
guesswork, starting from a certain analytic expression in
section~\ref{section.ricochet.proof} below. It is synthetically
constructed as follows: 
\begin{itemize} 
\item 
Start with arbitrary points $A,B,C,D$ on the conic. 
\item Let $V$ denote the intersection point of the tangents at
  $A$ and $C$, and let $F$ be on the conic such that $V,D,F$ are collinear. 
\item Let $W$ denote the intersection point of $AF$ and $CD$. 
\item Now mark off $Z$ on the conic such that $V, B, Z$ are collinear,
and finally $E$ such that $W,Z,E$ are collinear. 
\end{itemize} 

\begin{figure}
\includegraphics[width=7cm]{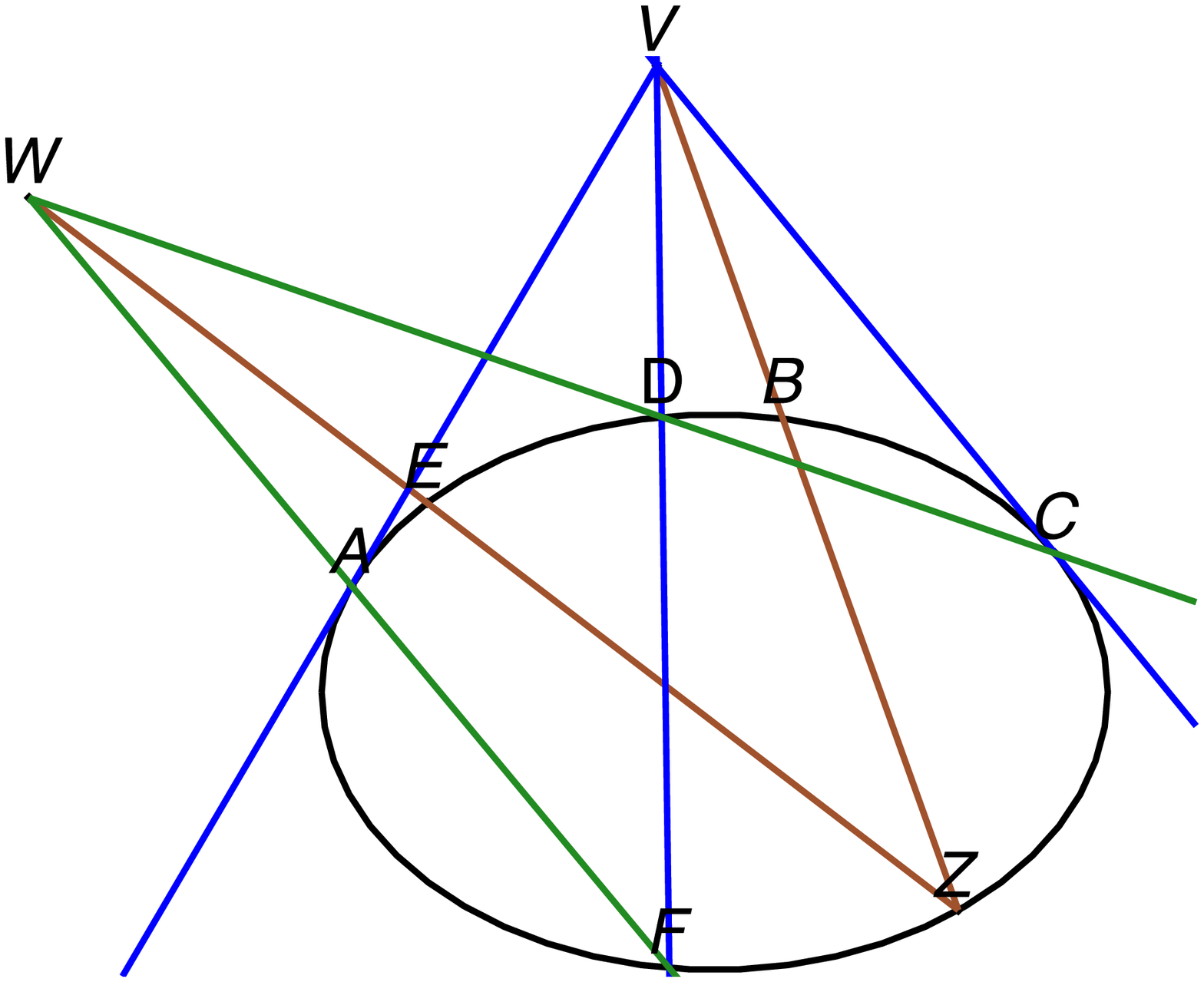}
\centerline{\small \textcolor{magenta}{\bf Diagram 3}}
\end{figure} 

In this situation, the Pascals 
\begin{equation} 
\pasc{A}{B}{C}{F}{E}{D}, \quad \pasc{A}{E}{C}{D}{B}{F} 
\label{ricochet.display} \end{equation} 
coincide; this will be proved in section~\ref{section.ricochet.proof}
below. (The common line is in fact $VW$, but the diagram would 
become too baroque for comprehension if any 
further lines were added to it.) One can imagine $B$ being struck by $V$ in the
direction of $Z$, bouncing off the conic and getting redirected to $E$, hence the term `ricochet'.

To recapitulate the main theorem, every sextuple of distinct points whose Pascals are not
all distinct must come from either Diagram 2 or Diagram 3. One 
can construct Diagram 2 starting from an arbitrary choice of $Q$ 
together with three lines through it, hence $\dim \InvH = 5$. 
Diagram 3 is completely determined by the
choice of $A, B, C, D$, hence $\dim \RicoV = 4$. 

The proof of the main theorem uses a case-by-case
analysis on pairs of Pascals, and each case is then disposed off using Gr{\"o}bner
basis computations. All such computations were carried out in {\sc Maple}. 

\subsection{} 
The next two sections are devoted to preliminaries. In
section~\ref{section.pascal.labels}, we recall the classical labelling schema for
Pascals. It is a beautiful combinatorial phenomenon which implicitly
involves the unique outer automorphism of the symmetric group on six objects. 

The group of automorphisms of $\P^2$ which preserve $\conic$ (not
necessarily pointwise, but as a set) is isomorphic to
$\text{PSL}(2,\complex)$. This group acts on all of the varieties mentioned above, and hence it is
convenient to use the language of binary forms and
$SL_2$-representations throughout (see section~\ref{section.sextics}). I have
included rather more explanation than what would have sufficed for
this paper alone, since I should like to refer to it in possible sequels to this paper. 

\medskip 

The literature on Pascal's theorem is very large. One of the 
best surveys of the field is due to George Salmon
(see~\cite[Notes]{SalmonConics}). The labelling schema, and a great deal of
other classical material is explained by H.~F.~Baker in his note `On the 
\emph{Hexagrammum Mysticum} of Pascal' in~\cite[Note II]{Baker}. An
engaging graphical presentation of this subject may be found at the {\sc url} 
\[ \text{\tt http://www.math.uregina.ca/{\raise.17ex\hbox{$\scriptstyle\sim$}}fisher/Norma/paper.html} \] 
maintained by J.~Chris Fisher and Norma Fuller. We
refer the reader to~\cite{KK, PedoeProj} for foundational notions in 
projective geometry, and to~\cite{Harris} for those in 
algebraic geometry.

\section{The Labelling Schema for Pascals} 
\label{section.pascal.labels} 

Start with the following sets 
\[ \six = \{1,2,3,4,5,6\}, \quad \text{and} \quad \ltr=
\{\bA,\bB,\bC,\bD,\bE,\bF\}. \] 
(The elements of $\ltr$ will eventually stand for points on
the conic, but at the moment they are pure letters.) 
A number duad is a $2$-element subset of $\six$, e.g.,
$\{3,5\}$. A number syntheme is a partition of $\six$ into three
number duads, e.g., $\{\{1,3\},\{2,6\},\{4,5\}\}$. We will flatten out the duads and synthemes for 
readability, i.e., write them as $35$ and $13.26.45$ etc. 
There are similar notions of a letter duad and a letter syntheme
answering to the set $\ltr$. For instance, $\bA\bE$ is a letter duad,
and $\bA\bC.\bD\bE.\bB\bF$ is a letter syntheme. 

Consider the sets $ND, NS, LD, LS$ of number duads, number synthemes,
letter duads, and letter synthemes respectively. Each of these four
sets has cardinality $15$. Now consider the following artfully
constructed diagonally symmetric table: 

\[ 
\begin{array}{|c|c|c|c|c|c|c|} \hline 
{} & \bA & \bB & \bC & \bD & \bE & \bF \\ \hline 
\bA & {} & 14.25.36 & 16.24.35 & 13.26.45 & 12.34.56 & 15.23.46 \\ 
\bB & 14.25.36 & {} & 15.26.34 & 12.35.46 & 16.23.45 & 13.24.56 \\ 
\bC & 16.24.35 & 15.26.34 & {} & 14.23.56 & 13.25.46 & 12.36.45 \\ 
\bD & 13.26.45 & 12.35.46 & 14.23.56 & {} & 15.24.36 & 16.25.34 \\ 
\bE & 12.34.56 & 16.23.45 & 13.25.46 & 15.24.36 & {} & 14.26.35 \\ 
\bF & 15.23.46 & 13.24.56 & 12.36.45 & 16.25.34 & 14.26.35 & {} \\ \hline 
\end{array} \] 

\bigskip 

A direct verification shows that it defines a bijection $LD \lra NS$;
where for instance, $\bB\bC$ is mapped to $15.26.34$. 

\subsection{} 
This table can be used to create a label for each Pascal. For
instance, consider the array $\arr{A}{E}{F}{C}{B}{D}$. Picture it as 

\begin{figure}[H]
\includegraphics[width=4cm]{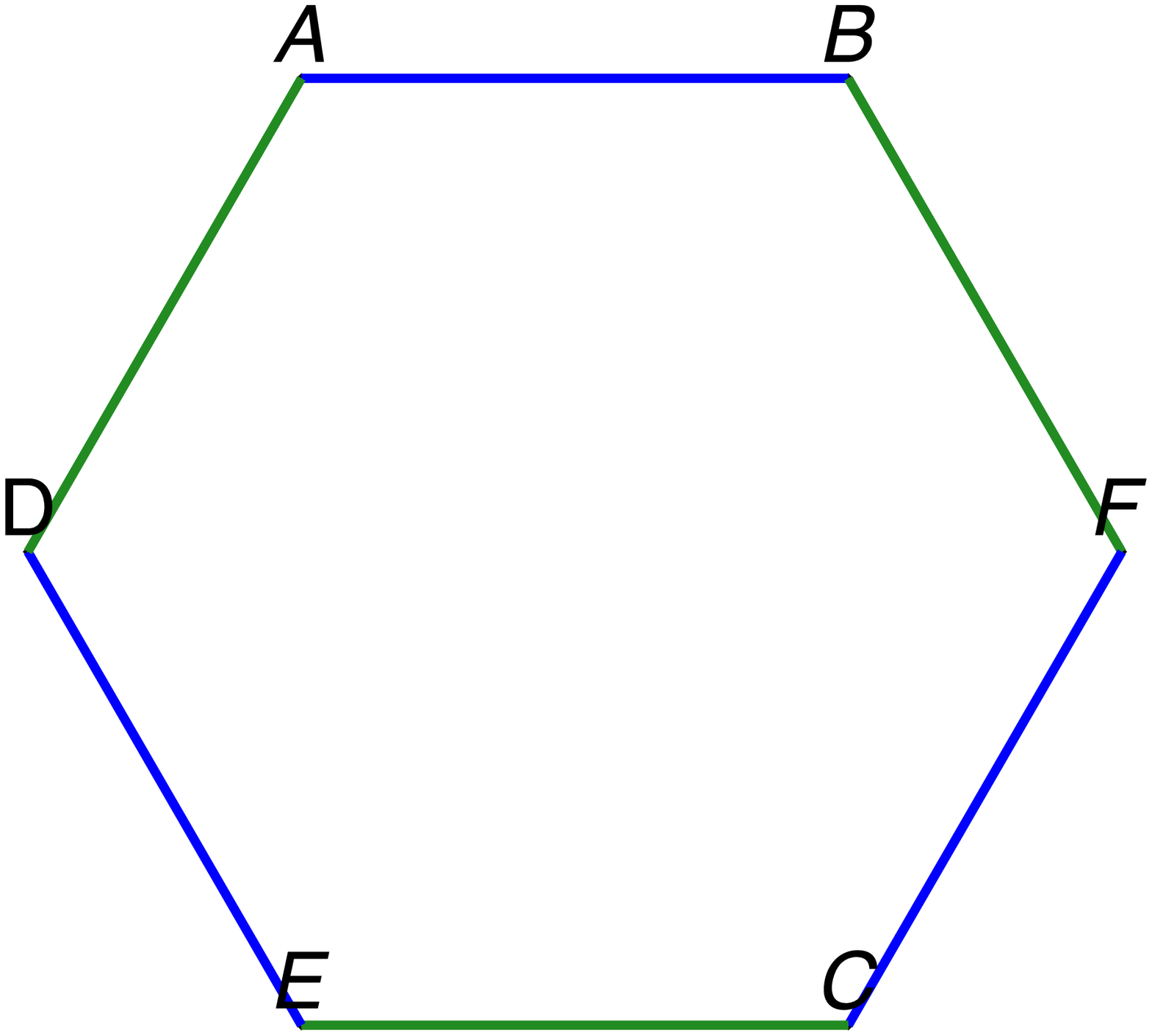}
\centerline{\small \textcolor{magenta}{\bf Diagram 4}}
\end{figure} 

so that each cross-hair intersection is between a blue and a green 
line forming opposite sides of the hexagon. 
Use the table above to find the number synthemes corresponding to the
blue lines: 
\[ \bA \bB \leadsto 14.25.36, \quad \bF \bC \leadsto 12.36.45, \quad 
\bE \bD \leadsto 15.24.36,  \] 
all of which have the duad $36$ in common. Similarly, those
corresponding to the green lines 
\[ \bA \bD \leadsto 13.26.45, \quad \bE \bC \leadsto 13.25.46, \quad 
\bF \bB \leadsto 13.24.56, \] 
have the duad $13$ in common. These two duads share the $3$, which
alternately combines with $1$ and $6$. Hence the corresponding Pascal 
$\pasc{A}{E}{F}{C}{B}{D}$ is given the label $k(3,16)$ or
$k(3,61)$. In summary, starting from an array of points, use the table
to extract two duads in the pattern $ab, ac$; and then the
corresponding Pascal is labelled $k(a,bc)$ or $k(a,cb)$. Since $a
\in \six$, and $\{b,c\} \subseteq \six \setminus \{a\}$, there are
altogether $6 \times \binom{5}{2} = 60$ labels, as they should be. 

The reader may wish to check that the Pascals in (\ref{inv.display})
are respectively
\[ k(1,23), \quad k(4,23), \quad k(5,23), \quad k(6,23). \] 
Those in (\ref{ricochet.display}) are respectively $k(1,23)$ and
$k(1,45)$. 

\subsection{} 
In the reverse direction, say we are given the label $k(2,35)$.  In
order to construct the corresponding array, start with the duads $23, 25$. Look for $23$ in
the table; it appears in positions $\bA\bF, \bB\bE,
\bC\bD$. Similarly, $25$ appears in $\bA\bB, \bC\bE,
\bD\bF$. This determines the hexagon: 

\begin{figure}[H]
\includegraphics[width=4cm]{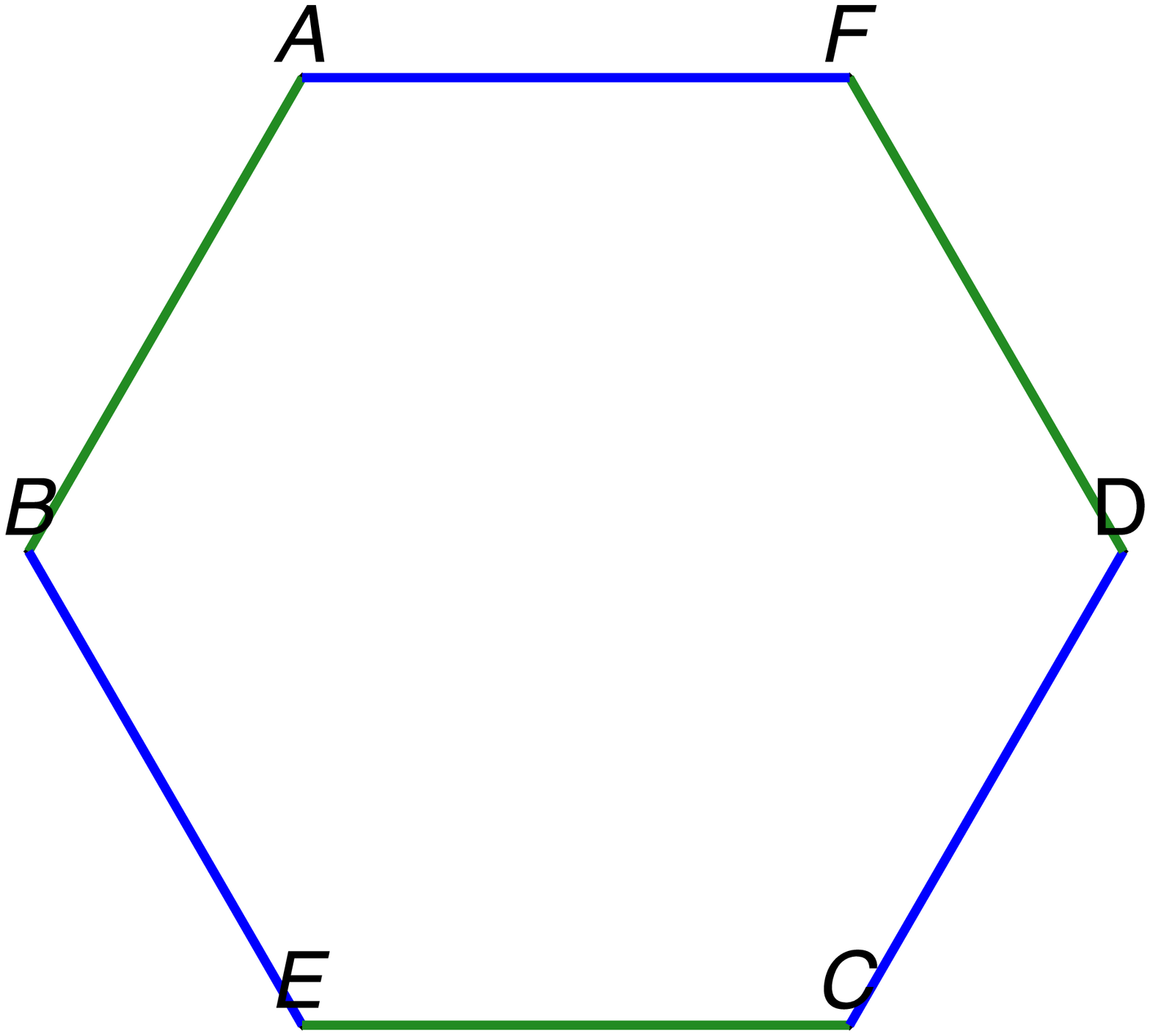}
\centerline{\small \textcolor{magenta}{\bf Diagram 5}}
\end{figure} 

and hence the array\footnote{It is of course understood that the hexagon is determined only
  up to rotation and reflection, and the array up to a permutation of
  rows and columns.} 
as $\arr{A}{D}{E}{C}{B}{F}$. In other words, the same table defines a bijection $ND \lra LS$, which takes 
$23$ to $\bA\bF.\bB\bE.\bC\bD$ etc., and then one can recover the
array from the images of the two duads. 
\subsection{} \label{section.map.SGs} 
Let $\SG(X)$ denote the symmetric group on the set $X$. 
Then the table defines an isomorphism $\SG(\ltr) \lra
\SG(\six)$. For instance, the image of the transposition $(\bA \, \bB)$ is
the product $(1 \, 4) \, (2 \, 5) \, (3 \, 6)$, and the map extends by
writing an arbitrary element as a product of transpositions. If we 
identify $\ltr$ and $\six$ as $\bA \leadsto 1,
\bB \leadsto 2, \dots, \bF \leadsto 6$, then this gives an 
outer automorphism $\omega$ of $\SG(\six)$, which is completely
specified\footnote{We follow the convention that the cycle $(1 \, 2 \,
  \dots 6)$ takes $1$ to $2$ etc.} by 
\[ (1 \, 2) \stackrel{\omega}{\lra} (1 \, 4) \, (2 \, 5) \, (3 \, 6), \qquad 
(1 \, 2 \, 3 \, 4 \, 5 \, 6)  \stackrel{\omega}{\lra} (2 \, 3 \, 6) \,
(4 \, 5). \] 
(Note that it does not preserve the cycle structure, and hence cannot
be inner.) A theorem of H{\"o}lder
characterises the outer automorphism groups of all finite symmetric
groups (see~\cite[Ch.~7]{Rotman}); it says that 
\[ \text{Out}(\SG(\{1,2,\dots,d\})) \simeq \begin{cases} 
{\mathbf Z}_2 & \text{if $d = 2,6$}, \\ 
\{e\}  & \text{otherwise.} 
\end{cases} \] 
Thus, $\omega$ represents the unique nontrivial element in
$\text{Out}(\SG(\six))$. A different identification of $\ltr$ with
$\six$ would amount to composing $\omega$ with an inner automorphism. 

The table above (along with its heavily Greek terminology of duads and
synthemes) was in essence constructed by Sylvester
(see~\cite{Sylvester}); however, I did not find his papers easy to
follow. What is usually called the \emph{Hexagrammum
Mysticum} is a much richer configuration than merely the Pascal lines,
and includes the Kirkman points and Cayley-Salmon lines
etc. They can all be labelled using the same formalism, and their
incidence relations can be read off from the labelling -- see
the note by Baker referred to above.  Other geometric perspectives on the outer automorphism may be
found in~\cite{Howard_etal}. 

\section{Binary Forms and Involutions} 
\label{section.sextics} 

In this section we will recast the necessary geometric notions in the language of binary
forms and $SL_2$-representations. A similar set-up is used
in~\cite{ego1}, where rather more detailed explanations are given. 

\subsection{} 
Let $V$ denote a two-dimensional complex vector space with basis $\ux
= \{x_1, x_2\}$, and a natural action of the group $SL(V)$. For $m \ge 0$, let $S_m$ denote
the $(m+1)$-dimensional vector space of homogeneous order $m$ forms in
$\ux$. It is an irreducible representation of $SL(V)$. 
Given integers $m, n \ge 0$ and $0 \le r \le \min(m,n)$, we have transvectant morphisms 
\[ S_m \otimes S_n \lra S_{m+n-2r}, \quad U \otimes
V \lra (U,V)_r;  \] 
given by the explicit formula 
\begin{equation} (U,V)_r = \frac{(m-r)! \, (n-r)!}{m! \, n!} \, 
\sum\limits_{i=0}^r \, (-1)^i \binom{r}{i} \, 
\frac{\partial^r U}{\partial x_1^{r-i} \, \partial x_2^i} \, 
\frac{\partial^r V}{\partial x_1^i \, \partial x_2^{r-i}} \, . 
\label{trans.formula} \end{equation} 

There is a symbolic calculus for transvectants, which is thoroughly
explained in~\cite[Ch.~1]{GrYo}. The basic theory of
$SL_2$-representations may be found in~\cite[Ch.~11]{FultonHarris}. 

\subsection{} 
Throughout, we will work inside the projective plane $\P S_2 \simeq \P^2$;
thus a nonzero quadratic form $Q \in S_2$ represents a point $[Q] \in
\P^2$. Its polar line is defined to be 
\[ \ell_Q = \{ [R] \in \P S_2: (R,Q)_2=0 \}. \] 
Every line in $\P^2$ is the polar of a unique point,
called its pole. There is a canonical isomorphism of $\P S_2$ with the dual plane $(\P
S_2)^*$, which maps $[Q]$ to $\ell_Q$. 

Given $Q, R \in S_2$, we have $(R,Q)_2 = (Q,R)_2$. Hence $[R] \in \ell_Q$ iff
$[Q] \in \ell_R$. The line of intersection of $[Q]$ and $[R]$ is given
by the polar of $[(Q,R)_1]$, and the point of intersection of $\ell_Q$
and $\ell_R$ is $[(Q,R)_1]$. 

\subsection{} 
Consider the Veronese imbedding 
\begin{equation} 
\P S_1 \stackrel{\phi}{\lra} \P S_2, \qquad 
[u] \lra [u^2]. \label{veronese} \end{equation} 
The image of $\phi$ is a smooth conic
$\conic$. If $Q = a_0 \, x_1^2 + a_1 \, x_1 \, x_2 + a_2 \, x_2^2$, then 
\[ (Q,Q)_2 = - \frac{1}{2} \, (a_1^2 - 4 \, a_0 \, a_2). \] 
Hence, 
\[ [Q] \in \conic \iff \text{$Q$ is the square of a linear form} \iff 
(Q,Q)_2=0 \iff [Q] \in \ell_Q. \] 

If $Q \in S_2$ factors as $u_1 \, u_2$, then the
points of intersection of $\ell_Q$ with $\conic$ are $\phi(u_1),
\phi(u_2)$. Dually, the tangent to the conic at either
$\phi(u_i)$ passes through $[Q]$. 

\subsection{}
A sextuple of unordered points $\Gamma = \{\phi(u_1), \dots, 
\phi(u_6)\}$ on $\conic$ will correspond to the binary sextic form 
$G_\Gamma = \prod\limits_{i=1}^6 u_i$, distinguished up to a scalar. Alternately, a nonzero form $G$
in $S_6$ will give a sextuple $\Gamma_G$
on $\conic$. This gives an isomorphism of $\P S_6$ with
$\Sym^6(\conic)$, where the discriminant hypersurface $\Delta
\subset\P S_6$ corresponds to sextuples with repeated points. It will
be occasionally convenient to use affine co-ordinates on $\conic$, by identifying 
$\phi(x_1 - \alpha \, x_2)$ with $\alpha$, and $\phi(x_2)$ with $\infty$. 

Since all incidences and intersections in $\P^2$ can be expressed as
transvectants, Pascal's theorem itself can be seen as a
transvectant identity (see~\cite[Theorem 2]{Leitenberger}). 
Define a \emph{hexad} to be an injective map $\ltr \stackrel{h}{\lra}
\conic$. We will write $h(\bA) = A, \dots, h(\bF) = F$, for the
corresponding points on $\conic$. If $\Hex$ denotes the set of all
hexads and $\Label_k$ the set of all labels, then we have a morphism 
\[ \Hex \lra \prod\limits_{\Label_k} \;  (\P^2)^*, \] 
which maps the hexad to its Pascals. The groups
$\SG(\ltr), \SG(\six)$ respectively act on $\Hex$ and the direct
product compatibly via the isomorphism in section~\ref{section.map.SGs}. Passing to quotients
by these actions, we get a morphism 
\[ \P^6 \setminus \Delta \lra \Sym^{60} \, (\P^2)^*, \] 
which maps a sextuple to the set of its Pascals. 
For what it is worth, I have calculated all the
Pascals for the sextuple $\Gamma = \{0,1,\infty,3,-5,7\}$ using {\sc Maple}, and verified that
they are in fact distinct. Hence, they must remain so for a
general $\Gamma$. 

\subsection{The quadratic involution} 
Fix a point\footnote{Henceforth we write $Q$ for $[Q]$ etc.~when no
  confusion is likely.}  $Q \in \P S_2$ away from $\conic$. It defines
an order $2$ automorphism (i.e., an involution) $\sigma_Q$ of $\conic$ as follows: 
if $\sz \in \conic$, then $\sigma_Q(\sz)$ is the other point of
intersection of $\conic$ with the line $Q\sz$. Now $\sigma_Q^2(\sz) =
\sz$, and $\sigma_Q(\sz) = \sz$ exactly when $Q\sz$ is tangent to
$\conic$. 
If $u \in S_1$ is such that $\phi(u) = \sz$, then $\sigma_Q(\sz)$
corresponds to the linear form $(Q,u)_1$. All of this is pursued further
in~\cite{ACquad}. 

Now $\sigma_Q$ extends to an involution of $\P^2$ by the following recipe:
given $R \in \P^2$, let $\sz_1, \sz_2$ be the (possibly coincident)
points where the polar of $R$ intersects $\conic$. Then define 
$\sigma_Q(R)$ to be the pole of the line joining $\sigma_Q(\sz_1)$ and
$\sigma_Q(\sz_2)$. In terms of transvectants, 
\[ \sigma_Q(R) = (Q,Q)_2 \, R - 2 \, (Q,R)_2 \, Q. \] 
Since $\sigma_Q(R)$ is a linear combination of $Q$ and $R$,
the points $Q, R, \sigma_Q(R)$ are collinear. The set of fixed points
of $\sigma_Q$ is $Q$ itself, together with the polar line of $Q$. 
(Thus, $\sigma_Q$ is a homology in the sense of~\cite[Ch.~11]{KK}). 

\subsection{} \label{section.inv}
Now assume that we have a hexad $\{A, \dots, F\}$ such that 
\[ \sigma_Q(A) = F, \quad \sigma_Q(B) = E, \quad \sigma_Q(C) = D,  \] 
as in Diagram 2. Consider the Pascal $\pasc{A}{B}{C}{F}{E}{D}$. Since
$\sigma_Q$ interchanges the lines $AE$ and $BF$, it must leave their
intersection point invariant. Similarly, $\sigma_Q$ leaves
each of the cross-hair intersections invariant, and hence they must
all lie on the polar of $Q$. It makes no difference to the argument if we select any one column in
the array and interchange its entries. We have proved the following proposition. 
\begin{Proposition} \sl 
With notation as above, each of the Pascals 
\[ \pasc{A}{B}{C}{F}{E}{D}, \quad \pasc{A}{B}{D}{F}{E}{C}, \quad 
\pasc{F}{B}{C}{A}{E}{D}, \quad \pasc{A}{E}{C}{F}{B}{D} \] 
is equal to the polar line of $Q$. 
\label{prop.inv.coincidence} \end{Proposition} 

As mentioned earlier, these Pascals carry labels $k(r,23)$ for $r \in
\{1,4,5,6\}$. By renaming the points, one would in general obtain four lines in the pattern 
\[ k(r,ab), \quad r \in \six \setminus \{a,b\}. \]

\subsection{The involutive hypersurface} 
A sextuple of points $\Gamma = \{ \sz_1, \dots, \sz_6 \}$ is said to be in involution if it is left invariant
by $\sigma_Q$ for some $Q \in \P^2$, and then $Q$ is said to be its
centre of involution. (In other words, the sextuple should fit into Diagram
2 for some $Q$.) Consider the variety 
\[ \InvH = \{ [G] \in \P^6 \setminus \Delta: \text{$\Gamma_G$ is in
  involution} \}. \] 

Change variables so that $Q = x_1 \,x_2$. If $\sz \in \conic$
corresponds to $u = x_1 + \alpha \, x_2$, then $\sigma_Q(\sz)$
corresponds to\footnote{Henceforth we will write $\Box$ for a
  multiplicative scalar whose precise value is unimportant. For
  instance, $\Box$ stands for $- \frac{1}{2}$ here.} 
$(Q,u)_1 = \Box \, (x_1 - \alpha \, x_2)$, and then $u \, (Q, u)_1$ is
a quadratic with no $x_1 \, x_2$ term. Thus $\Gamma_G$
is in involution with respect to $Q$, if and only if $G$ can be
written as a form in $x_1^2, x_2^2$. In other words, $\InvH$ is the
variety of sextic forms which can be written as 
\begin{equation} c_1 \, u_1^6 + c_2 \, u_1^4 \, u_2^2 + c_3 \, u_1^2
  \, u_2^4 + c_4 \, u_2^6, \qquad (c_i \in \complex), \label{inv.form} \end{equation}
for some linear forms $u_1, u_2$ (cf.~\cite[\S 260]{SalmonAlg}). 

\subsection{The covariants of a binary sextic} 
\label{section.sextic.covariants} 
The complete minimal system of covariants of a generic binary sextic is given in
\cite[p.~156]{GrYo}. We will not reproduce it here; but only note down a
few of its members which are relevant to the subject at hand. 

Let $G$ denote a generic sextic, and write $\vartheta_{m,q}$ for
a covariant of degree-order $(m,q)$. This means that,
when written out in full, 
\[ \vartheta_{m,q} = \sum\limits_{i=0}^q \, \theta_i \, x_1^{q-i} \,
x_2^i, \] 
where $\theta_i$ are homogeneous forms of degree $m$ in the coefficients of
$G$. If $q=0$, then $\vartheta_{m,0}$ is called an invariant of degree
$m$. Now define 
\begin{equation} 
\vartheta_{2,4} = (G, G)_4, \quad 
\vartheta_{3,2} = (G, \vartheta_{2,4})_4, \quad 
\vartheta_{8,2} = (\vartheta_{2,4},\vartheta_{3,2}^2)_3, \quad 
\vartheta_{15,0} = ((G, \vartheta_{2,4})_1,\vartheta_{3,2}^4)_8.
\label{covariants.sextic} \end{equation} 

It is known that $\InvH$ is a hypersurface defined by the vanishing of 
$\vartheta_{15,0}$ (see~\cite[\S 4.10]{ACquad}). Moreover, 
$\vartheta_{8,2}$ evaluated on the form~(\ref{inv.form}) gives 
$\Box \, u_1 \, u_2$, which is $Q$. Thus, if $G$ is in involution, then $\vartheta_{8,2}$ can be used to
`detect' its centre if it is unique. (However, if $G$ is arbitrary,
then $\vartheta_{8,2}$ has no geometric meaning that I know of.) 
As we will see in section~\ref{section.triple.symmetry},
it may happen that a sextuple in a highly special position has more
than one centre of involution, and then $\vartheta_{8,2}$ vanishes
identically. 

I have programmed the transvectant formula~(\ref{trans.formula}) in
{\sc Maple}, so that these covariants can be calculated on a specific
$G$ wherever necessary. 

\subsection{The ricochet configuration} \label{section.ricochet.proof}
Assume that the hexad $\{A, \dots, F\} \subseteq \conic$ is in ricochet
configuration as shown in Diagram 3.  

\begin{Proposition} \sl 
Both the Pascals $\pasc{A}{B}{C}{F}{E}{D}, \pasc{A}{E}{C}{D}{B}{F}$
coincide with the line $VW$. 
\label{prop.ricochet.coincidence.proof} 
\end{Proposition} 
\proof 
This is a straightforward computation with transvectants. Choose co-ordinates such that 
\[ A = \phi(x_1), \quad C = \phi(x_2), \quad B = \phi(x_1 - x_2),
\quad D = \phi(x_1 - d \, x_2). \] 
Then $V = \Box \, x_1 \, x_2$, and $F$ corresponds to $(V,x_1 - d \, x_2)_1 =
\Box \, (x_1 + d \, x_2)$. Hence 
\[ W = (x_1 \, (x_1 + d \, x_2), x_2 \, (x_1 - d \, x_2))_1 = \Box \, (x_1^2 -
2 \, d \, x_1 \, x_2 - d^2 \, x_2^2). \] 
Now $Z$ is given by $(x_1 \, x_2, x_1 - x_2)_1 = \Box \, (x_1 + x_2)$,
and finally $E$ by 
\[ (W,x_1 + x_2)_1 = \Box \, (x_1 + \frac{d^2-d}{d+1} \,
x_2). \] 
One can similarly calculate all the cross-hair intersections and the
lines joining them. It turns out that either Pascal is given by the quadratic
form $P = x_1^2 + d^2\, x_2^2$; or in other words, it is the polar of
$[P]$. Since $(P,V)_2 = (P,W)_2=0$, it must pass through $V$ and
$W$. \qed 

\medskip 

Notice that $P$ factors as $(x_1 + d \, x_2 \, \sqrt{-1} ) \, (x_1 -
d \, x_2 \, \sqrt{-1})$, i.e., if $VW \cap \conic = \{I,J\}$, then 
$I,J$ have affine co-ordinates $\pm \, d \, \sqrt{-1}$. This implies that we
have cross-ratios 
\[ \langle \, A,C, I, J \, \rangle = \langle \, D, F, I, J \, \rangle = -1, \] 
i.e., $I,J$ is a harmonically conjugate pair with respect to $A,C$ as
well as $D,F$. Since $V,W$ are determined by $A,C,D$, the
common Pascal is independent of the position of $B$. These
observations suggest that a more conceptual and less computational proof of this proposition should be
possible, but I do not see one. 

\section{The Main Theorem} 
In this section we will establish the following theorem. 
\begin{Theorem} \sl 
Let $\Gamma$ be a hexad, and assume that 
$s,t$ are two labels such that $k(s) = k(t)$ for $\Gamma$. Then $\Gamma$ is
either in involution or in ricochet configuration. 
\label{main.theorem} \end{Theorem} 
\proof
After applying an automorphism of $\conic$, we may assume that the points of
$\Gamma$ are given in affine co-ordinates as 
\begin{equation} A = 0, \quad B = 1, \quad C = \infty, \quad D = p, \quad E = q, \quad
F = r, 
\label{affine.sextuple} \end{equation} 
and hence 
\[ G_\Gamma = x_1 \, (x_1 - x_2) \, x_2 \, (x_1 - p \, x_2) \, (x_1 -
q \, x_2) \, (x_1 - r \, x_2). \] 

Now the proof simply goes through all possible $s$ and $t$, but one 
can introduce a small technical device to reduce the number of cases. 

\subsection{} 
Given a label $s = (a, bc)$, write $s' = \{a\}$, and $s'' =
\{b,c\}$. For two labels $s, t$, define their interference matrix 
\[ I_{st} = \left[ \begin{array}{cc} s' \cdot t' & s' \cdot t'' \\ s'' \cdot
   t' &
s'' \cdot t'' \end{array} \right], \] 
where $s' \cdot t''$ means the cardinality of the set $s' \cap
t''$ and so on. 

For instance, if $s = (1,23), t = (2, 36)$, then 
\[ s' = \{1\}, \quad s'' = \{2, 3\}, \quad t' = \{2\}, \quad t'' =
\{3, 6\}, \quad \text{and} \quad I_{st} = \left[ \begin{array}{cc} 0 & 0 \\ 1 & 1 \end{array}
\right]. \] 

After applying a permutation of $\six$, we may assume once and for all that
$s = (1,23)$. It corresponds to the array $\arr{A}{B}{C}{F}{E}{D}$, and then a direct
calculation as in section \ref{section.ricochet.proof}  shows that $k(1,23)$ is given by the 
quadratic form 
\begin{equation}
 (q-r) \, x_1^2 + (p \, r-p \, q + p-q) \, x_1 \, x_2 + r \, (q-p)
\, x_2^2. 
\label{q.k123} \end{equation}
 If $t,u$ are two labels such that $I_{st} = I_{su}$, then one can find a
permutation carrying $t$ into $u$ which preserves $s$, hence it
suffices to consider any one example of $t$ for any given interference
matrix. The following are all the possibilities for $I_{st}$. 

\[ \begin{array}{lllll} 
I^{(1)} = \intarr{1}{0}{0}{0}, & I^{(2)} = \intarr{1}{0}{0}{1}, &
I^{(3)} = \intarr{0}{0}{1}{0}, & I^{(4)} = \intarr{0}{0}{1}{1}, & I^{(5)} = \intarr{0}{1}{1}{1}, \\
I^{(6)} = \intarr{0}{1}{1}{0}, & I^{(7)} = \intarr{0}{0}{0}{0}, &
I^{(8)} = \intarr{0}{0}{0}{1}, & 
I^{(9)} = \intarr{0}{0}{0}{2}. \end{array} \] 

Since the whole question is symmetric in $s$ and $t$, it is
unnecessary to consider the transpose of $I^{(3)}$ or $I^{(4)}$.

\subsection{} 
Let $I_{st} = I^{(2)} = \intarr{1}{0}{0}{1}$. We may assume
$t=(1,24)$, corresponding to the array $\arr{A}{D}{F}{C}{E}{B}$. A very
similar calculation shows that $k(t)$ is given by
\begin{equation} (p-r) \, x_1^2 + (r- p \, q) \, x_1 \, x_2 + p \, r \, (q - 1) \,
x_2^2. \label{q.k124} \end{equation}
If $k(s)=k(t)$, then (\ref{q.k123}) and (\ref{q.k124}) must be
scalar multiples of each other, and hence the $2 \times 3$ matrix of their
coefficients must have all of its minors zero. This gives a system of
polynomial equations in $p,q,r$. One solves it by finding a Gr{\"o}bner
basis of the resulting ideal, after imposing an elimination
order on the variables (see~\cite[Ch.~2]{AL} or~\cite[Ch.~3]{CLO} for the
technique). However,  in this case, the only solutions are 
\[ p=r=0, \quad q=1, r=0, \quad q=p, r=0, \quad p=q=1, \quad q=1,
r=p, \quad p=q=r.\] 
None of these is legal, since each would force $\Gamma$ to have
a repeated point. We conclude that the two Pascals cannot
coincide. Similarly, we get no legal solutions for $I^{(j)}, j = 3, 6, 7,8$. 

\subsection{} The remaining four cases are geometrically more
interesting. They have the common feature that apart from illegal
solutions as above (which will not be explicitly mentioned), there is
a unique nontrivial solution in every case. 

Say $I_{st} = I^{(4)} = \intarr{0}{0}{1}{1}$, then we may take $t = (2,
34)$ corresponding to the array $\arr{A}{B}{D}{E}{C}{F}$. A similar 
calculation gives the solution 
\[ q = \frac{p}{p+1}, \quad r = \frac{p}{1-p^2}, \] 
with $p$ arbitrary. (It is, of course, subject to the constraint that 
no two points of $\Gamma$ should coincide, which excludes only
finitely many values of $p$. Henceforth this proviso is tacitly understood
whenever we have free parameters.) Substitute the solution into $G =
G_\Gamma$ to get a binary sextic whose coefficients are functions of
$p$. Now a rather long calculation using the formulae in
(\ref{covariants.sextic}) shows that
$\vartheta_{15,0}(G)=0$, hence $\Gamma$ must be in involution. The
centre of the involution is found to be 
\[ \vartheta_{8,2}(G) = Q = \Box \, (x_1^2 - 2 \, p \, x_1 \, x_2 +
\frac{p^2}{1+p} \, x_2^2). \] 
The lines $AE, CD, BF$ pass through $Q$. Hence, by Proposition \ref{prop.inv.coincidence}, the Pascals 
\begin{equation}
 \pasc{A}{B}{C}{E}{F}{D}, \quad 
\pasc{A}{B}{D}{E}{F}{C}, \quad 
\pasc{A}{F}{C}{E}{B}{D}, \quad 
\pasc{E}{B}{C}{A}{F}{D} 
\label{pascals234} \end{equation} 
all coincide with each other; or what is the same, 
$k(4,56) = k(1,56) = k(2,56) = k(3,56)$. Thus we have the curious
situation that if $k(1,23), k(2,34)$ coincide, then four other Pascals
are also forced to coincide. 

Here is a more geometric way to see this configuration: fix $Q,
A,B,E,F$, and allow the line $CD$ to pivot around $Q$. 

\begin{figure}[H]
\includegraphics[width=4cm]{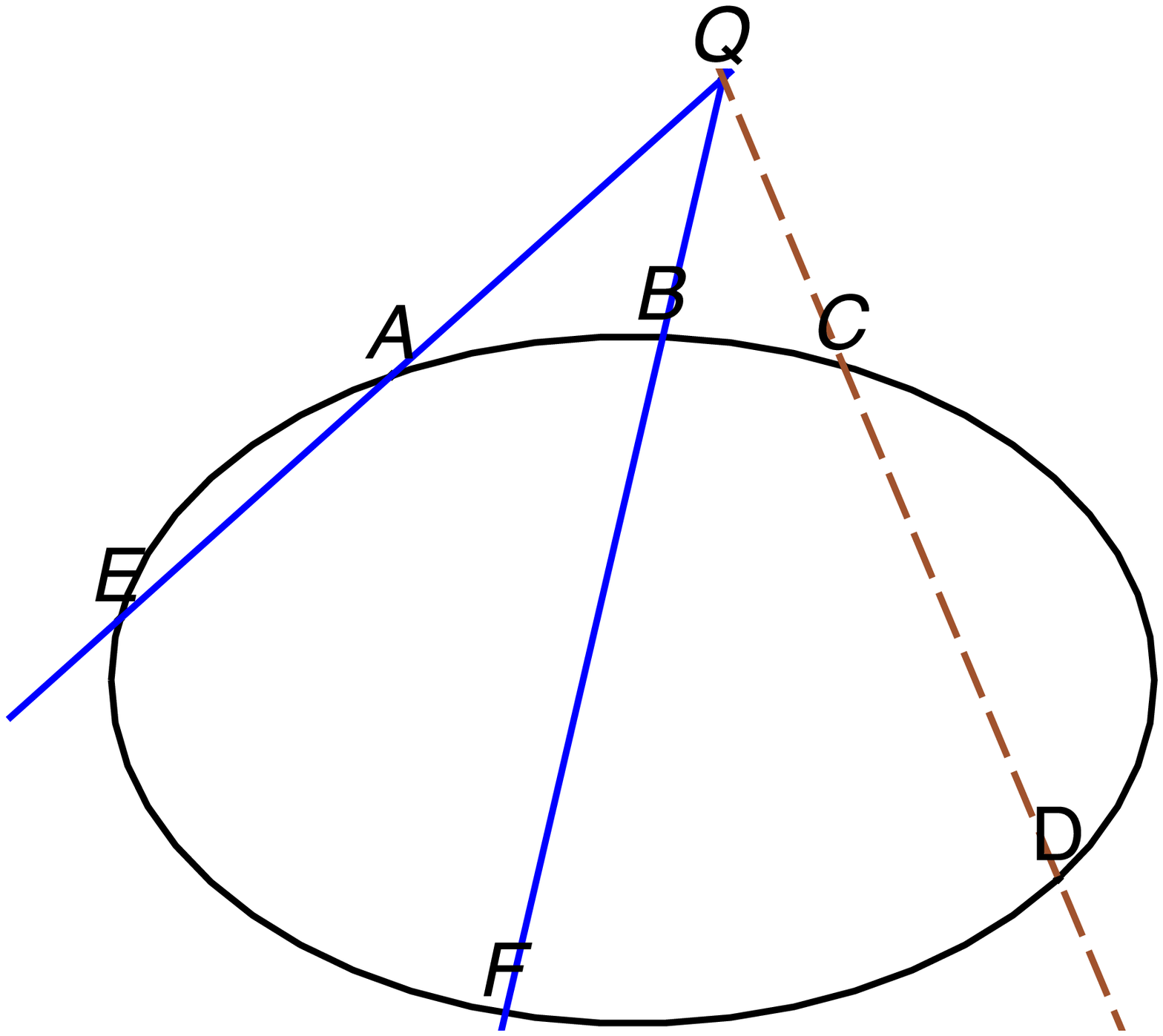}
\centerline{\small \textcolor{magenta}{\bf Diagram 6}}
\end{figure}  

The Pascals in~(\ref{pascals234}) coincide for any position of
$CD$. Furthermore, 
\[ k(1,23) \leadsto \underbrace{\pasc{A}{B}{C}{F}{E}{D}}_{\lambda_1},
\qquad k(1,24) \leadsto \underbrace{\pasc{A}{B}{D}{E}{C}{F}}_{\lambda_2} \] 
both pass through $Q=AE \cap BF = BF \cap CD$.  Let $\Pi_Q$ denote the
pencil of lines through $Q$; then we have a two-to-one morphism 
\[ \conic \stackrel{g_1}{\lra} \Pi_Q, \qquad C \lra \lambda_1 \] 
which maps $C$ to the line joining $BD \cap CE$ with $Q$. The 
similar morphism 
\[ \conic \stackrel{g_2}{\lra} \Pi_Q, \qquad C \lra \lambda_2 \] 
maps $C$ to the line joining $AC \cap BE$ with $Q$. Since 
$\Pi_Q \simeq \P^1$ has a unique rational double cover up to
isomorphism\footnote{This may be seen as follows: such a cover is
  completely determined by its two simple branch points, and any two
  points on $\P^1$ can be taken to any other by the Fundamental
Theorem of Projective Geometry.}, there must be an automorphism $\tau$ of $\Pi_Q$
such that $\tau \circ g_1 = g_2$. But then $\tau$ must have at least
one fixed point (in fact generically two such points), that is to
say, a line $\lambda \in \Pi_Q$ such that $\tau(\lambda) =
\lambda$. Hence, fixed points of $\tau$ correspond to positions of $C$
such that $\lambda_1 = \lambda_2$.

\subsection{} Assume that $I_{st} = I^{(9)} = \intarr{0}{0}{0}{2}$,
then we may take $t=(4,23)$. Using the procedure above, one gets the
two parameter solution 
\[ q = \frac{p \, (r-1)}{p-1}, \] 
with $p,r$ arbitrary. Then one finds that $\vartheta_{15,0}(G)=0$, and 
$\vartheta_{8,2}(G) = Q = x_1^2 - 2 \, p \, x_1 \, x_2 + p \, r
\, x_2^2$. A calculation shows that $AF, BE, CD$ intersect in $Q$, and we are
simply in the generic involutive configuration of section~\ref{section.inv}. 

\subsection{} \label{section.triple.symmetry} 
Assume that $I_{st} = I^{(5)} = \intarr{0}{1}{1}{1}$, then we may take
$t=(2,13)$. The same procedure gives the one-parameter solution 
\[ q = \frac{p-1}{p}, \quad r = \frac{1}{1-p}. \] 
Now $\vartheta_{15,0}(G)=0$, hence $\Gamma$ must be in
involution. However, $\vartheta_{8,2}(G)$ also vanishes identically, hence
one should look for multiple centres. On the other hand, substituting
the solution into~(\ref{q.k123}) shows that $k(1,23)$ is given by 
\[ T = x_1^2 - x_1 \, x_2 + x_2^2 = (x_1 + \theta \, x_2) \, (x_1 +
\theta^2 \, x_2), \qquad \theta = e^{\frac{2 \pi \sqrt{-1}}{3}} \] which is
independent of $p$. The factors of $T$ are suggestive of a connection with
`equi-anharmonicity', i.e., the phenomenon where the cross-ratio of
four points on a line admits a threefold symmetry
(see~\cite[Ch.~II.8]{Yaglom}). Indeed, it turns out that the cyclic
group ${\mathbf Z}_3$ acts on the entire structure in such a way that, 
four distinct groups of Pascals coincide amongst themselves. 

Consider the linear transformation $\sigma$ of $S_1$ which acts by 
\[ x_1 \lra x_1-x_2, \quad x_2 \lra x_1. \] 
It induces an action on $\P S_2$ and $\conic$, either of which will
also be denoted by $\sigma$. Notice that $\sigma^3$ is the scalar multiplication
by $-1$, and hence acts as the identity on $\P S_2$. It is easy to
check that the action of $\sigma$ on $\conic$ stabilizes the set $\Gamma = \{A, \dots, F \}$,
and acts as the permutation $(A \, B \, C) \, (D \, F \, E)$. (That is
to say, $\sigma$ takes $A$ to $B$, and $D$ to $F$ etc.) 

Define points 
\[ M = \phi(x_1 + \theta \, x_2), \quad N = \phi(x_1 + \theta^2 \,
x_2), \] 
on $\conic$, then $\sigma(M) = M, \sigma(N) = N$, and hence
the line $MN$ (which is the polar of $T$) is fixed (as a set) by
$\sigma$. Note the cross-ratios
\[ \langle \, C,A,B,M \, \rangle = \langle \, \infty, 0, 1, - \theta
\, \rangle = - \theta, \qquad 
\langle \, C, A, B, N \, \rangle = \langle \, \infty,0,1,-\theta^2
\rangle = - \theta^2;  \] 
which agrees with the fact that 
\[ \langle C, A, B, M \, \rangle = \langle \, \sigma(C), \sigma(A),
\sigma(B), \sigma(M) \, \rangle = \langle \, A, B, C, M \, \rangle,  \] 
and similarly for $N$. In classical terminology, $\{C,A,B,M\}$ and
$\{C,A,B,N\}$ are equi-anharmonic tetrads. 

Now let $\alpha = p-1, \beta = 1, \gamma = -p$, and consider the three quadratic
forms: 
\[ 
Q_6 = \alpha \, x_1^2 + 2 \, \beta \, x_1 \, x_2 + \gamma \, x_2^2,
\quad 
Q_4 = \beta \, x_1^2 + 2 \, \gamma \, x_1 \, x_2 + \alpha \,
x_2^2, \quad 
Q_5 = \gamma \, x_1^2 + 2 \, \alpha \, x_1 \, x_2 + \beta \,
x_2^2. \] 
(Notice the cyclic movement of $\alpha, \beta, \gamma$.) Then 
$(Q_6,T)_2 = (Q_4,T)_2 = (Q_5,T)_2 = 0$, 
and hence all $[Q_i]$ are on the line $MN$. The action of $\sigma$ on $\P^2$ is such that 
$[Q_6] \ra [Q_4] \ra [Q_5] \ra [Q_6]$. 
A simple check shows that the lines $AD, BE, CF$ intersect in $Q_6$; furthermore $AE,CD,
BF$ intersect in $Q_4$, and $AF,CE,BD$ in $Q_5$. Thus $\Gamma$ is a highly
special configuration which is in involution with respect to three
different centres. 

\begin{figure}[H]
\includegraphics[width=10cm]{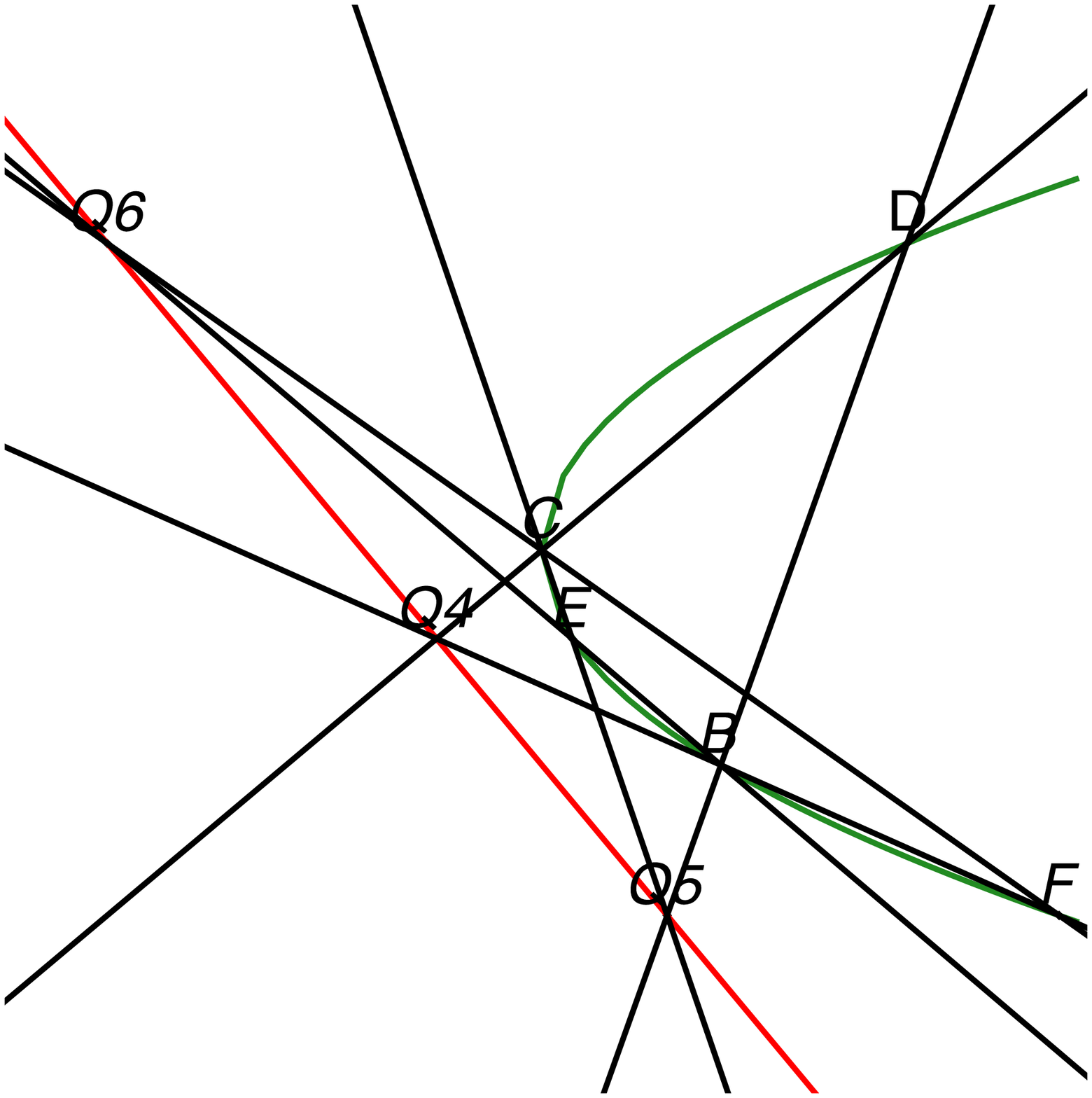}
\centerline{\small \textcolor{magenta}{\bf Diagram 7}}
\smallskip 
\parbox{12cm}{\tiny The point $A$ (not shown) is to the far right at
 infinity. The points $M,N$, not being real, cannot be shown.} 
\end{figure}  
 
By Proposition~\ref{prop.inv.coincidence}, we have
the following sets of coincidences: 
\begin{equation} \begin{aligned} 
{} & k(1,45) = k(2,45)=k(3,45)=k(6,45), \\
& k(1,56) = k(2,56)=k(3,56)=k(4,56), \\ 
& k(1,46) = k(2,46)=k(3,46)=k(5,46). \end{aligned} 
\label{three.pascals} \end{equation} 
Or, what comes to the same thing, the map $\SG(\ltr) \lra
\SG(\six)$ sends $(\bA \, \bB \, \bC) \, (\bD \, \bF \, \bE)$ to $(4
\, 5 \, 6)$; the latter induces a cyclic action on the three groups of Pascals
in~(\ref{three.pascals}), and also explains the subscripts in $Q_i$. 

We are yet to explain the identity $k(1,23) = k(2,13)$. Notice that
$k(1,23) \leadsto \pasc{A}{B}{C}{F}{E}{D}$ must pass
through $AD \cap CF = Q_6$. Applying $\sigma$ to the points, 
\[ \pasc{A}{B}{C}{F}{E}{D} \stackrel{\sigma}{\lra} 
\pasc{B}{C}{A}{E}{D}{F} = \pasc{A}{B}{C}{F}{E}{D}, \] 
that is to say, $k(1,23)$ is left invariant by $\sigma$. However, it
must pass through $\sigma(Q_6) = Q_4$, and hence must be the line $Q_6Q_4
= MN$. By the same argument, either of the Pascals 
\[ k(2,13) \leadsto \pasc{A}{B}{C}{D}{F}{E}, \qquad 
k(3,12) \leadsto \pasc{A}{B}{C}{E}{D}{F} \] 
is also equal to $MN$, and thus $k(1,23)=k(2,13)=k(3,12)$.

\subsection{} \label{section.richochet.proof} 
There remains the case $I_{st} = I^{(1)} =
\intarr{1}{0}{0}{0}$. Assuming $t = (1,45)$, we get the solution 
\begin{equation} q = p \, (1-p)/(1+p), \quad r = -p, 
\label{analytic.soln.ricochet} \end{equation}
with $p$ arbitrary. 

It turns out that $\vartheta_{15,0}(G)$ does not vanish as a function
of $p$, hence $\Gamma$ is not in involution for generic $p$. (However, see
section~\ref{section.intersect.RY} below.) But 
notice that if we substitute this analytic solution 
into~(\ref{affine.sextuple}), everything agrees exactly with the proof of
Proposition~\ref{prop.ricochet.coincidence.proof}, with $p$ in place of
$d$. This shows that $\Gamma$ is in ricochet configuration, hence 
the proof of Theorem~\ref{main.theorem} is complete. \qed  

As mentioned earlier, I used (\ref{analytic.soln.ricochet}) as a starting point, and only
afterwards reached the construction in
section~\ref{section.ricochet.construction}. Several false steps were necessary before it was
found. 

It would be interesting to have an essentially synthetic proof of the main
theorem, i.e., one which uses as much classical projective geometry and as little
explicit calculation as possible. 

\subsection{} 
Given an interference pattern $I$, one may consider the variety 
\[ \Omega_I = \{ [G] \in \P S_6 \setminus \Delta: 
\text{The sextuple $\Gamma_G$ has coincident Pascals in pattern
  $I$} \}. \] 
These are $SL_2$-equivariant subvarieties of $\P^6 \setminus \Delta$,
and it would be of interest to find their degrees, desingularisations,
and defining equations. As we have
seen, $\Omega_{I^{(j)}}$ is empty for $j=2,3,6,7,8$, and $\Omega_{I^{(9)}} =
\InvH$. In any of the remaining cases we get a one-parameter solution
in $p$, and since the $SL_2$-orbit of $\Gamma$ for a specific $p$ is
three-dimensional (see~\cite{AF}), the variety $\Omega_I$ itself
must be four-dimensional. It is contained in $\InvH$ for $j=4,5$, but not for
$j=1$. 

I tried to calculate the ideal of the `ricochet locus' $\RicoV = \Omega_{I^{(1)}}$
inside the co-ordinate ring of $\P^6$ using elimination of 
variables (rather as in~\cite[\S 4.8]{ACquad}), but could not get
the computation to terminate. This is unfortunately a chronic difficulty
with practical elimination theory. Even so, a direct calculation
with the fundamental system of sextics shows that there
is one invariant in degree $6$, and two independent invariants in
degree $10$ vanishing on this locus. One can at least
conclude that the ideal is not a complete intersection. 

\subsection{} \label{section.intersect.RY} 
The value of the invariant $\vartheta_{15,0}$ on the `ricochet' form is: 
\[ p^{18} \, (p^2+3) \, (3 \, p^2+1) \, (p^2 + 1) \, (p^2+p+1) \,
(p^2-p+1) \, (p^2+2 \, p-1)^2 \, (p^2- 2 \, p-1)^2 (p-1)^3 (p+1)^3. \] 
It vanishes for finitely many $p$, hence the intersection $\RicoV
\cap \InvH$ is a finite union of $SL_2$-orbits. 

\section{Pascals on the Discriminant Locus} 

Hitherto we have assumed that $\Gamma$ consists of six distinct points,
but all the Pascals are well-defined if any one pair of points is allowed
to come together. 

\subsection{} \label{section.double.six} 
In order to see this, assume that $A=B$, and $C, D, E, F$ are distinct from each other and
from $A$. We will interpret $AB$ as the tangent to $\conic$ at
$A$. Given an array of points, one may assume that $A$ occupies the top
left corner, and then it is only necessary to consider the following three positions of $B$. 
\begin{equation} 
\underbrace{\arr{A}{B}{D}{F}{E}{C}}_{I}, \quad \underbrace{\arr{A}{C}{D}{B}{F}{E}}_{II}, \quad 
\underbrace{\arr{A}{C}{D}{F}{B}{E}}_{III}. \label{three.types} \end{equation}

In case I, $AE \cap BF = A$ and the other two cross-hair
intersections are on the line $AC$, hence the Pascal is $AC$. 

In case II, $AF \cap BC, AE \cap BD$ both equal $A$,
hence the Pascal is the line joining $A$ to $CE \cap DF$. 

In order to see that the Pascal is well-defined in case III, it
is enough to show that the points $P = AB \cap CF, P' = AE \cap DF$
cannot coincide. If they did, $AP$ would be tangent
to the conic at $A$ and would contain $E$,  which is impossible. 

\subsection{} 
However, if $\Gamma$ has either a threefold point or two double
points, then some of the Pascals become undefined. If $A=B=C$,
then $\pasc{A}{B}{C}{F}{E}{D}$ is no longer defined, since all
cross-hair intersections are at $A$. If $A=B$ and $C=D$,
then $\pasc{A}{B}{E}{C}{D}{F}$ becomes undefined, since the line $AC =
AD \cap BC$ may not contain the point $AF \cap CE$. 

\subsection{} 
If $\Gamma \in \Delta$, then it is already clear
that many of the Pascals must coincide; for instance, in case I above,
the Pascal remains the same for all permutations of 
$D,E,F$. In this section we will describe all such coincidences. 

The general  picture is that the set of
labels splits into three types I, II, III as in
(\ref{three.types}). Type I splits further into $4$ classes with $6$
elements each, type II into $3$ classes with $4$ elements each, and
type III into $12$ classes with $2$ elements each. Altogether there are $19$
equivalence classes, such that all Pascals in each class are equal. 
For a general $\Gamma$ in $\Delta$, these $19$ lines are distinct. 

{\bf Type I:} All Pascals of the form $\pasc{A}{B}{\star}{\star}{\star}{C}$ are
equal, which gives a $6$-element equivalence class. To determine their
labels, note that we know two of the sides of the corresponding
hexagon, namely $AC, BC$. From the table, 
\[ \bA\bC \leadsto 16.24.35, \quad \bB\bC \leadsto 15.26.34. \] 
The label must come from two duads (i.e., one from each number syntheme)
which have an element in common. The pair $16, 15$ leads to $k(1,56)$,
and similarly the other possibilities are 
\[ k(6,12), \quad k(2,46), \quad k(4,23), \quad k(5,13), \quad k(3,45). \] 
We get three similar equivalence classes by replacing $C$ with
$D, E, F$. 

{\bf Type II:} Consider all arrays of the form 
$\arr{A}{\star}{\star}{B}{\star}{\star}$, where the rightmost $2 \times
2$ block is one of 
\[ \left[ \begin{array}{cc} C & D \\ F & E \end{array} \right], \quad 
\left[ \begin{array}{cc} C & F \\ D & E \end{array} \right], \quad 
\left[ \begin{array}{cc} D & E \\ C & F \end{array} \right], \quad 
\left[ \begin{array}{cc} F & E \\ C & D \end{array} \right]. \] 
The Pascal is the line joining $A$ to $CE \cap DF$ in all cases,
hence we have a $4$-element equivalence class. The labels are easily
determined to be $k(4,36), k(1,36), k(3,14), k(6,14)$. They are
constructed on the following model: start with two number duads $ab, cd$ having
no element in common (here $14, 36$), and then combine them as 
\[ k(a,cd), \quad k(b,cd), \quad k(c,ab), \quad k(d,ab). \] 
We get two more such classes from $CD \cap EF$ and $CF \cap
DE$. Since $\bA\bB \leadsto 14.25.36$, picking any two duads out of the
three will give one of the three equivalence classes. 

{\bf Type III:} We have 
\[ \pasc{A}{C}{D}{F}{B}{E} = \pasc{B}{C}{D}{F}{A}{E},  \] 
or what is the same, $k(2,15) = k(5,24)$. The latter Pascal may 
be written as $\pasc{A}{F}{E}{C}{B}{D}$, hence in
general we have a $2$-element equivalence class consisting of 
\[ \pasc{A}{P_1}{Q_1}{P_2}{B}{Q_2} \quad \text{and} \quad 
\pasc{A}{P_2}{Q_2}{P_1}{B}{Q_1}, \] 
where $\{ P_1, P_2, Q_1, Q_2 \} = \{C,D,E,F\}$. There are
$\frac{4!}{2} = 12$ such classes. Their labels are formed on the
following pattern: from the image of $\bA\bB \leadsto 14.25.36$, pick any of
the three duads (say $ab$), pick another (say $cd$) and now form the
$2$-element class of $k(a,bc), k(b,ad)$. (Note that the construction
is not symmetric in $ab, cd$, nor in $c,d$.) 

\subsection{} 
In order to assert that there are no further coincidences for a general
$\Gamma$ in $\Delta$, it is sufficient to check this on one example. After
choosing, 
\[ A = B = 0, \quad C = \infty, \quad D = 1, \quad E = -2, \quad F
=3, \] 
I have calculated all the Pascals, and verified that there are
precisely $19$ of them. 

In conclusion, if $\T$ denotes the locus of 
sextic forms which have at least a triple root or two double roots,
then we have a morphism $\P^6 \setminus \T \lra \Sym^{60} \,
(\P^2)^*$ just as in section~\ref{section.pascal.map}. By the main theorem, the
preimage of the big diagonal is $\Delta \cup \InvH
\cup \RicoV$. According to standard procedure, one can 
blow up $\P^6$ along $\T$ to extend the morphism
(see~\cite[Ch.~II.7]{Hartshorne}); but I will leave this analysis for a sequel.

\medskip 

\centerline{--} 

\vspace{1cm}

\parbox{7cm}{ \small 
Jaydeep Chipalkatti \\
Department of Mathematics \\ 
University of Manitoba \\ 
Winnipeg, MB R3T 2N2 \\ 
Canada. \\ \\
{\tt chipalka@cc.umanitoba.ca}}

\end{document}